\documentclass[12pt]{article}
\usepackage{latexsym}
\usepackage{amsmath}
\usepackage{amssymb}
\usepackage{amsfonts}
\usepackage[colorlinks=true]{hyperref}
\hypersetup{urlcolor=blue, citecolor=red}

\numberwithin{equation}{section}

\newtheorem{theorem}{Theorem}[section]

\newcommand{\beq}{\begin{equation}}
\newcommand{\eeq}{\end{equation}}
\def\nm{\noalign{\medskip}}

\newcommand{\eps}{\varepsilon}
\newcommand{\R}{{\mathbb R}}

\renewcommand{\div}{{\rm div}}

\newcommand{\la}{\langle}
\newcommand{\ra}{\rangle}
\newcommand{\BR}{\bf R}
\newcommand{\Bn}{\bf n}
\newcommand{\Bt}{\bf t}

\title{Bounds on the volume fraction of the two-phase shallow shell using one measurement}
\author{Hyeonbae Kang\thanks{Department of Mathematics, Inha University, Incheon 402-751, Korea (hbkang@inha.ac.kr)} \and
Graeme W. Milton\thanks{Department of Mathematics, University of Utah, Salt Lake City, UT 84112, USA (milton@math.utah.edu)}
\and
Jenn-Nan Wang\thanks{Department of Mathematics, NCTS (Taipei), National Taiwan University, Taipei 106, Taiwan (jnwang@math.ntu.edu.tw)}}

\date{}

\begin{document}
\maketitle

\begin{abstract}
We study the size estimate problem for the two-phase shallow shell equations in this paper. Our aim is to derive bounds on the volume fraction of each phase assuming that the material properties of the two phases are given. The approach in this paper is based on the translation method. One of the key steps is to connect the shallow shell equations to the thin plate equation. 
\end{abstract}


\section{Introduction}\label{sec1}
\setcounter{equation}{0}

Let an elastic shallow shell consist of two materials. We can assume that one of them is an inclusion whose size is unknown. A fundamental question is to estimate the volume fraction of the inclusion by a noninvasive method such as measuring boundary responses. In this paper we want to prove some bounds on the volume fraction based on the translation method.  To begin, we let $\Omega$ be a bounded simply connected domain in $\R^2$ with smooth boundary $\partial\Omega$.  For a shallow shell, its middle surface is described by
$\{(x_1,x_2,\eps{\theta}(x_1,x_2)):(x_1,x_2)\in\overline{\Omega}\}$
for $\eps>0$, where $\eps$ is a small parameter.  Let ${\bf u}=(u_1,u_2,u_3)=({\bf
u}',u_3):\Omega\to\R^3$ represent the displacement vector of the
middle surface. Taking $\eps\to 0$ with appropriate scalings, Ciarlet and Miara \cite{cm} showed that ${\bf u}$ satisfies the following equations:
\begin{equation}\label{shell}
\begin{cases}
\div\ {\bf s}^{\theta}=0\quad\text{in}\quad\Omega,\\
\div\div\ {\bf m}-\div({\bf s}^{\theta}\nabla\theta)=0\quad\text{in}\quad\Omega,
\end{cases}
\end{equation}
where the constitutive laws are given as
\begin{eqnarray}\label{form}
{\bf m}&=&\frac{4\mu}{3}\nabla^2u_3+\frac{4\lambda\mu}{3(\lambda+2\mu)}(\text{Tr}\nabla^2u_3)I_2,\notag\\
{\bf s}^{\theta}&=&4\mu
{\bf e}^{\theta}+\frac{4\lambda\mu}{\lambda+2\mu}(\text{Tr}{\bf e}^{\theta})I_2,\\
{\bf e}^{\theta}&=&\frac
12\left (\nabla{\bf u}'+(\nabla{\bf u}')^T+\nabla\theta\otimes\nabla u_3+(\nabla\theta\otimes\nabla u_3)^T\right),\notag
\end{eqnarray}
and $\lambda,\mu$ are Lam\'e coefficients. Here, ${\bf a}\otimes{\bf b}$ denotes the tensor product of vectors ${\bf a}$ and ${\bf b}$ and $I_2$ is the $2\times 2$ identity matrix.  The Cauchy data on $\partial\Omega$ associated with \eqref{shell} is given by
\begin{equation}\label{data1}
\begin{cases}
({\bf u}';{u}_3, u_{3,n})\ (\text{Dirichlet}),\\
\left({\bf s}^{\theta}{\bf n};\
(\div\ {\bf m}-{\bf s}^{\theta}\nabla\theta)\cdot{\bf n}+({\bf m}{\bf n}\cdot{\bf t})_{,t}, {\bf m}{\bf n}\cdot{\bf n}\right)\ (\text{Neumann}),
\end{cases}
\end{equation}
where ${\bf n}$ is the boundary normal vector field and ${\bf t}$ is the tangential vector field in the positive orientation. Also, we define $g_{,n}=\nabla g\cdot{\bf n}$ and $g_{,t}=\nabla g\cdot{\bf t}$ for any function $g$. Here and elsewhere in the paper ${\bf C}{\bf a}\cdot {\bf b}$ means $({\bf C}{\bf a})\cdot {\bf b}$ for any $2\times 2$ matrix ${\bf C}$. It is clear that $(\div\ {\bf m}-{\bf s}^{\theta}\nabla\theta)\cdot{\bf n}+({\bf m}{\bf n}\cdot{\bf t})_{,t}$ can be replaced by $(\div\ {\bf m})\cdot{\bf n}+({\bf m}{\bf n}\cdot{\bf t})_{,t}$ since ${\bf s}^{\theta}{\bf n}$ is already given.

Assume now that the shell consists of two phases, i.e.,
\[
\lambda=\lambda_1\chi_1+\lambda_2\chi_2\quad\text{and}\quad\mu=\mu_1\chi_1+\mu_2\mu_2,
\]
where
\[
\chi_j=\begin{cases}1\quad\text{in phase}\ j,\\ 0\quad\text{otherwise}.\end{cases}
\]
Here we are interested in the inverse problem of estimating the volume fraction of the first phase $f_1$ (or the second phase $f_2$) by one pair of Cauchy data \eqref{data1} when $\lambda_1,\lambda_2,\mu_1,\mu_2$ are given.  The approach taken here is based on the variational method, precisely, the translation method. On one hand,  this paper extends results in \cite{kkm} (for conductivity) and in \cite{ml} (for elasticity), which used the translation method of Murat and Tartar \cite{mt1,mt2,mt3} and Lurie and Cherkaev \cite{lc1,lc2}, to the shallow shell equation. On the other hand, we would like to give another application of the authors' recent result on the equivalence of the inverse boundary problems for 2D elasticity and the thin plate \cite{kmw}.  Using the method in \cite{kmw}, the idea here is to transform the elasticity-like equation of \eqref{shell} into a plate-like equation (4th order). At the same time, the Cauchy data for the new equation are completely determined by the original one. Therefore, the original problem is reduced to the size estimate problem for a decoupled system of 4th order equations, which is easier to deal with. For the shallow shell equations with an
isotropic inhomogeneous medium in which the medium inside of the inclusion is not known, bounds on the volume fraction involving constants which are not explicitly given were derived in \cite{dlw}. The method used in \cite{dlw} is based on some quantitative uniqueness estimates for the shallow shell equations.

This paper is organized as follows. In Section~2, we transform the first equation of \eqref{shell} into a 4th order equation and discuss the corresponding Cauchy data for the new system. The inverse problem is then investigated in Section~3.

\section{Transformation of the equation and the Cauchy data correspondence}

It was known that the 2D elasticity equation and the thin plate equation are equivalent (see for example \cite[Section~2.3]{mbook}). The transformation of the first equation of \eqref{shell} into a plate-like equation will follow the same idea. Denote
\[
{\bf s}^{\theta}=\begin{pmatrix}s_{11}^{\theta}&s_{12}^{\theta}\\s_{12}^{\theta}&s_{22}^{\theta}\end{pmatrix}.
\]
The first equation of \eqref{shell} is written as
\[
\begin{cases}
s_{11,1}^{\theta}+s_{12,2}^{\theta}=0,\\
s_{12,1}^{\theta}+s_{22,2}^{\theta}=0.
\end{cases}
\]
Hence, there exists a function $\psi$ such that
\begin{equation}\label{11}
{\bf
s}^{\theta}=\begin{pmatrix}s_{11}^{\theta}&s_{12}^{\theta}\\s_{21}^{\theta}&s_{22}^{\theta}\end{pmatrix}=\begin{pmatrix}\psi_{,22}&-\psi_{,12}\\-\psi_{,12}&\psi_{,11}\end{pmatrix}={\R}\nabla^2\psi,
\end{equation}
where the fourth order tensor ${\mathbb R}$ is defined by
\[
{\mathbb R}{\bf M} = {\bf R}_{\perp}^T {\bf M} {\bf R}_{\perp}
\]
for any $2\times 2$ matrix ${\bf M}$, where
\[
{\bf R}_{\perp}=\begin{pmatrix}0&1\\-1&0\end{pmatrix}.
\]

The second equation of \eqref{form} can be written as
\begin{equation}\label{12}
{\bf e}^{\theta}={\mathbb S}{\bf s}^{\theta},
\end{equation}
where the fourth order tensor ${\mathbb S}$ is defined by
\begin{equation}\label{14}
{\mathbb S}{\bf A}=\frac{1}{4\mu}{\bf
A}^{\text{sym}}-\frac{\lambda}{4\mu(3\lambda+2\mu)}(\text{Tr}{\bf
A})I_2
\end{equation}
for any matrix ${\bf A}$. On the other hand, from the form of ${\bf
e}^{\theta}$ (the third equation of \eqref{form}), we have
\begin{equation}\label{111}
{\bf e}^{\theta}-\frac{1}{2}\left(\nabla\theta\otimes\nabla u_3+(\nabla\theta\otimes\nabla u_3)^T\right)=\frac 12\left (\nabla{\bf u}'+(\nabla{\bf u}')^T\right)
\end{equation}
and it follows that
\begin{equation}\label{15}
\div\div\left({\mathbb R}{\bf e}^{\theta}\right)-\frac
12\div\div\left({\mathbb R}(\nabla\theta\otimes\nabla u_3)+{\mathbb
R}(\nabla\theta\otimes\nabla u_3)^T\right)=0.
\end{equation}
Using \eqref{11} and \eqref{12} in \eqref{15} yields
\begin{equation}\label{16}
\div\div({\mathbb L}\nabla^2\psi)-\frac 12\div\div\left({\mathbb
R}(\nabla\theta\otimes\nabla u_3)+{\mathbb
R}(\nabla\theta\otimes\nabla u_3)^T\right)=0,
\end{equation}
where ${\mathbb L}={\mathbb R}{\mathbb S}{\mathbb R}$. It is easily seen that ${\mathbb L}={\mathbb S}$. Replacing the
first equation of \eqref{shell} by \eqref{16}, \eqref{shell} is
transformed into
\begin{equation}\label{shell2}
\begin{cases}
\div\div({\mathbb L}\nabla^2\psi)-\frac 12\div\div\left({\mathbb
R}(\nabla\theta\otimes\nabla u_3)+{\mathbb
R}(\nabla\theta\otimes\nabla u_3)^T\right)=0,\\
\div\div({\mathbb M}\nabla^2u_3)-\div({\mathbb
R}\nabla^2\psi\nabla\theta)=0,
\end{cases}
\end{equation}
where $${\mathbb M}\nabla^2u_3={\bf m}=\frac{4\mu}{3}\nabla^2u_3+\frac{4\lambda\mu}{3(\lambda+2\mu)}(\text{Tr}\nabla^2u_3)I_2.$$

We now discuss the Cauchy data corresponding to \eqref{shell2}. The Cauchy data for \eqref{shell2} are given by
\begin{equation}\label{cd1}
\left\{\begin{split}
&(\psi,\psi_{,n});\\
&\div({\mathbb L}\nabla^2\psi)\cdot{\bf n}-\frac
12\div\left({\mathbb R}(\nabla\theta\otimes\nabla u_3)+{\mathbb
R}(\nabla\theta\otimes\nabla u_3)^T\right)\cdot{\bf n}\\
&+\left(({\mathbb L}\nabla^2\psi){\bf n}\cdot{\bf t}-\frac
12\left({\mathbb R}(\nabla\theta\otimes\nabla u_3)+{\mathbb
R}(\nabla\theta\otimes\nabla u_3)^T\right){\bf n}\cdot{\bf
t}\right)_{,t};\\
&({\mathbb L}\nabla^2\psi){\bf n}\cdot{\bf n}-\frac 12\left({\mathbb
R}(\nabla\theta\otimes\nabla u_3)+{\mathbb
R}(\nabla\theta\otimes\nabla u_3)^T\right){\bf n}\cdot{\bf n}
\end{split}\right.
\end{equation}
and
\begin{equation}\label{cd2}
\left\{\begin{split}
&(u_3,u_{3,n});\\
&\left(\div({\mathbb M}\nabla^2u_3)-{\mathbb R}\nabla^2\psi\nabla\theta\right)\cdot{\bf n}+(({\mathbb M}\nabla^2u_3){\bf n}\cdot{\bf t})_{,t};\\
&({\mathbb M}\nabla^2u_3){\bf n}\cdot{\bf n}.
\end{split}\right.
\end{equation}
The Cauchy data \eqref{cd2} correspond to the second equation of \eqref{shell2}, which are already given in \eqref{data1}. Now we want to investigate \eqref{cd1}. In view of authors' recent result \cite{kmw}, we can show that \eqref{cd1} are completely determined by $({\bf u}',{\bf s}^{\theta}{\bf n})$ on $\partial\Omega$. In fact, by \eqref{11}, one can determine $\psi$ and $\psi_{,n}$ from ${\bf s}^{\theta}{\bf n}$. On the other hand, let $\eps=\frac{1}{2} (\nabla{\bf u}' + (\nabla{\bf u}')^T)$, then we have
 \[
 {\mathbb R}\eps={\BR}_{\perp}^T \eps {\BR}_{\perp} = 
 \begin{pmatrix} u_{2,2} & -\frac{1}{2} u_{1,2} - \frac{1}{2} u_{2,1} \\
 -\frac{1}{2} u_{1,2} - \frac{1}{2} u_{2,1} & u_{1,1} \end{pmatrix}. 
 \]
Thus we obtain
 \begin{align*}
 \div \, {\BR}_{\perp}^T \eps {\BR}_{\perp} & =
 \begin{pmatrix} u_{2,12} -\frac{1}{2} u_{1,22} - \frac{1}{2} u_{2,12} \\
 \nm
 u_{1,12} -\frac{1}{2} u_{1,21} - \frac{1}{2} u_{2,11} \end{pmatrix}  
 = \begin{pmatrix} \frac{1}{2} u_{2,12} -\frac{1}{2} u_{1,22} \\
 \nm
 \frac{1}{2} u_{1,12} - \frac{1}{2} u_{2,11} \end{pmatrix} \\
 &= {\BR}_{\perp}^T \begin{pmatrix} (\frac{1}{2} u_{1,2} -\frac{1}{2} u_{2,1})_{,1} \\
 \nm
 (\frac{1}{2} u_{1,2} - \frac{1}{2} u_{2,1})_{,2} \end{pmatrix} 
 = {\BR}_{\perp}^T \nabla (\frac{1}{2} u_{1,2} - \frac{1}{2} u_{2,1}),
 \end{align*}
which implies
 \begin{align*}
 (\div \, {\BR}_{\perp}^T \eps {\BR}_{\perp}) \cdot{\Bn}
 &= ({\BR}_{\perp} {\Bn}) \cdot \nabla (\frac{1}{2} u_{1,2} - \frac{1}{2} u_{2,1}) \\
 &= - {\Bt} \cdot \nabla (\frac{1}{2} u_{1,2} - \frac{1}{2} u_{2,1}) \\
 &= (\frac{1}{2} u_{2,1} - \frac{1}{2} u_{1,2} )_{,t}.
 \end{align*}
Combining \eqref{11}, \eqref{12}, and \eqref{111} yields that
 \begin{equation}\label{2.8}
 \begin{split}
&\div({\mathbb L}\nabla^2\psi)\cdot{\bf n}-\frac
12\div\left({\mathbb R}(\nabla\theta\otimes\nabla u_3)+{\mathbb
R}(\nabla\theta\otimes\nabla u_3)^T\right)\cdot{\bf n}\\
&+\left(({\mathbb L}\nabla^2\psi){\bf n}\cdot{\bf t}-\frac
12\left({\mathbb R}(\nabla\theta\otimes\nabla u_3)+{\mathbb
R}(\nabla\theta\otimes\nabla u_3)^T\right){\bf n}\cdot{\bf
t}\right)_{,t}\\
&= (\frac{1}{2} u_{2,1} - \frac{1}{2} u_{1,2} )_{,t}+ \left(({\BR}_{\perp}^T \eps {\BR}_{\perp} {\Bn}) \cdot {\Bt}\right)_{,t}.
\end{split}
\end{equation}
Observe that $\frac{1}{2} (u_{2,1} - u_{1,2})$ can be expressed as 
 $$
 \frac{1}{2} (u_{2,1} - u_{1,2}) = \left({\BR}_\perp^T \begin{pmatrix} 0 & -\frac{1}{2} (u_{2,1} - u_{1,2}) \\
 \frac{1}{2} (u_{2,1} - u_{1,2}) & 0 \end{pmatrix} {\BR}_\perp {\Bn} \right) \cdot {\Bt}.
 $$ 
We also have 
 \begin{equation*}
 ({\BR}_{\perp}^T \eps {\BR}_{\perp} {\Bn}) \cdot {\Bt} = \left({\BR}_\perp^T \begin{pmatrix} u_{1,1} & \frac{1}{2} u_{1,2} + \frac{1}{2} u_{2,1} \\
 \frac{1}{2} u_{1,2} + \frac{1}{2} u_{2,1} & u_{2,2} \end{pmatrix} {\BR}_\perp {\Bn} \right) \cdot {\Bt}. 
 \end{equation*}
From \eqref{2.8} we thus obtain
\begin{equation}\label{Mt}
\begin{split}
&\div({\mathbb L}\nabla^2\psi)\cdot{\bf n}-\frac
12\div\left({\mathbb R}(\nabla\theta\otimes\nabla u_3)+{\mathbb
R}(\nabla\theta\otimes\nabla u_3)^T\right)\cdot{\bf n}\\
&+\left(({\mathbb L}\nabla^2\psi){\bf n}\cdot{\bf t}-\frac
12\left({\mathbb R}(\nabla\theta\otimes\nabla u_3)+{\mathbb
R}(\nabla\theta\otimes\nabla u_3)^T\right){\bf n}\cdot{\bf
t}\right)_{,t}\\
&= \left( \left({\BR}_\perp^T \begin{pmatrix} u_{1,1} & u_{1,2}\\
 u_{2,1} & u_{2,2} \end{pmatrix} {\BR}_\perp {\Bn} \right) \cdot {\Bt}\right)_{,t} \\
&= - ({\Bt} \cdot {\BR}_\perp^T (\nabla{\bf u}') {\Bt})_{,t} = -( ({\BR}_\perp {\Bt}) \cdot (\nabla {\bf u}') {\Bt})_{,t} = - ({\Bn} \cdot (\nabla {\bf u}') {\Bt})_{,t}=-({\bf n}\cdot{\bf u}'_{,t})_{,t}.
 \end{split}
 \end{equation}
 We can also deduce that 
 \beq\label{Mn}
 \begin{split}
 &({\mathbb L}\nabla^2\psi){\bf n}\cdot{\bf n}-\frac 12\left({\mathbb
R}(\nabla\theta\otimes\nabla u_3)+{\mathbb
R}(\nabla\theta\otimes\nabla u_3)^T\right){\bf n}\cdot{\bf n}\\
 &= {\Bn} \cdot ({\BR}_{\perp}^T \eps {\BR}_{\perp}) {\Bn} = {\Bt} \cdot \eps {\Bt} = {\Bt} \cdot (\nabla{\bf u}') {\Bt}={\bf t}\cdot{\bf u}'_{,t}.
 \end{split}
 \eeq
Formulae \eqref{Mt} and \eqref{Mn} show that the Neumann data in \eqref{cd1} are completely determined by ${\bf u}'_{,t}$.

It is helpful to take a closer look at $\nabla u_3$ in the Neumann data of \eqref{cd1}. It is obvious that $\nabla u_3$ (on $\partial\Omega$) is
determined by $u_3$ and $u_{3,n}$.  In other words, the values of
$\left({\mathbb R}(\nabla\theta\otimes\nabla u_3)+{\mathbb
R}(\nabla\theta\otimes\nabla u_3)^T\right){\bf n}\cdot{\bf t} $ and
$\left({\mathbb R}(\nabla\theta\otimes\nabla u_3)+{\mathbb
R}(\nabla\theta\otimes\nabla u_3)^T\right){\bf n}\cdot{\bf n}$ are
given in terms of $u_3$ and $u_{3,n}$. It remains to simplify
$$\div\left({\mathbb R}(\nabla\theta\otimes\nabla u_3)+{\mathbb
R}(\nabla\theta\otimes\nabla u_3)^T\right)\cdot{\bf n}.$$ Notice that
\[
{\mathbb R}(\nabla\theta\otimes\nabla
u_3)=\begin{pmatrix}\theta_{,2}u_{3,2}&-\theta_{,2}u_{3,1}\\-\theta_{,1}u_{3,2}&\theta_{,1}u_{3,1}\end{pmatrix}
\]
and
\[
{\mathbb R}(\nabla\theta\otimes\nabla
u_3)^T=\begin{pmatrix}\theta_{,2}u_{3,2}&-\theta_{,1}u_{3,2}\\-\theta_{,2}u_{3,1}&\theta_{,1}u_{3,1}\end{pmatrix}.
\]
Thus, we have that
\begin{eqnarray*}
&&\div({\mathbb R}(\nabla\theta\otimes\nabla u_3))\cdot{\bf n}\\
&=&\begin{pmatrix}(\theta_{,2}u_{3,2})_{,1}-(\theta_{,1}u_{3,2})_{,2}\\-(\theta_{,2}u_{3,1})_{,1}+(\theta_{,1}u_{3,1})_2\end{pmatrix}\cdot{\bf n}\\
&=&\begin{pmatrix}\theta_{,21}u_{3,2}-\theta_{,12}u_{3,2}\\-\theta_{,21}u_{3,1}+\theta_{,12}u_{3,1}\end{pmatrix}\cdot{\bf n}+(\theta_{,2}u_{3,21}-\theta_{,1}u_{3,22})n_1+(-\theta_{,2}u_{3,11}+\theta_{,1}u_{3,12})n_2\\
&=&\begin{pmatrix}\theta_{,21}u_{3,2}-\theta_{,12}u_{3,2}\\-\theta_{,21}u_{3,1}+\theta_{,12}u_{3,1}\end{pmatrix}\cdot{\bf n}+\theta_{,2}\left((u_{3,1})_{,2}n_1-(u_{3,1})_{,1}n_2\right)+\theta_{,1}\left(-(u_{3,2})_{,2}n_1+(u_{3,2})_{,1}n_2\right)\\
&=&\begin{pmatrix}\theta_{,21}u_{3,2}-\theta_{,12}u_{3,2}\\-\theta_{,21}u_{3,1}+\theta_{,12}u_{3,1}\end{pmatrix}\cdot{\bf n}+\theta_{,2}\left((u_{3,1})_{,2}t_2+(u_{3,1})_{,1}t_1\right)+\theta_{,1}\left(-(u_{3,2})_{,2}t_2-(u_{3,2})_{,1}t_1\right)\\
&=&\begin{pmatrix}\theta_{,21}u_{3,2}-\theta_{,12}u_{3,2}\\-\theta_{,21}u_{3,1}+\theta_{,12}u_{3,1}\end{pmatrix}\cdot{\bf n}+\theta_{,2}(u_{3,1})_{,t}-\theta_{,1}(u_{3,2})_{,t},
\end{eqnarray*}
which is determined by $u_{3}$ and $u_{3,n}$, where we write ${\bf n}=(n_1,n_2)^T$. Likewise, we can compute that
\begin{eqnarray*}
&&\div({\mathbb R}(\nabla\theta\otimes\nabla u_3)^T)\cdot{\bf n}\\
&=&\begin{pmatrix}(\theta_{,2}u_{3,2})_{,1}-(\theta_{,2}u_{3,1})_{,2}\\-(\theta_{,1}u_{3,2})_{,1}+(\theta_{,1}u_{3,1})_2\end{pmatrix}\cdot{\bf n}\\
&=&\begin{pmatrix}\theta_{,21}u_{3,2}-\theta_{,22}u_{3,1}\\-\theta_{,11}u_{3,2}+\theta_{,12}u_{3,1}\end{pmatrix}\cdot{\bf n}+(\theta_{,2}u_{3,21}-\theta_{,2}u_{3,12})n_1+(-\theta_{,1}u_{3,21}+\theta_{,1}u_{3,12})n_2\\
&=&\begin{pmatrix}\theta_{,21}u_{3,2}-\theta_{,22}u_{3,1}\\-\theta_{,11}u_{3,2}+\theta_{,12}u_{3,1}\end{pmatrix}\cdot{\bf n},
\end{eqnarray*}
which is also determined by $u_{3}$ and $u_{3,n}$. Having these computations, we conclude that the boundary data $\{\div({\mathbb L}\nabla^2\psi)\cdot{\bf n}+({\mathbb L}\nabla^2\psi){\bf n}\cdot{\bf t})_{,t},({\mathbb L}\nabla^2\psi){\bf n}\cdot{\bf n}\}$ entirely depends on ${\bf u}'$ and $(u_3,u_{3,n})$ on $\partial\Omega$.

For the purpose of studying the size estimate problem, we want to show that $\la\text{det}\nabla^2\psi\ra$, $\la\nabla^2\psi\ra$, and $\la{\mathbb L}\nabla^2\psi\ra$ are null-Lagrangians. Throughout the paper, $\la f\ra$ denotes the average of the field $f$, i.e.,
\[
\la f\ra=\frac{1}{|\Omega|}\int_{\Omega}f.
\]
Using the integration by parts, we immediately obtain that
\begin{equation*}
\la\nabla^2\psi\ra=\frac{1}{|\Omega|}\int_{\Omega}\nabla^2\psi=\frac{1}{|\Omega|}\int_{\partial\Omega}{\bf
n}(\nabla \psi)^T=\frac{1}{|\Omega|}\int_{\partial\Omega}{\bf
n}(\psi_{,t}{\bf t}+\psi_{,n}{\bf n})^T.
\end{equation*}
Hence $\la\nabla^2\psi\ra$ is a null-Lagrangian. Multiplying the first equation of \eqref{shell2} by $x_1^2/2, x_1x_2/2, x_2^2/2$, respectively, and using integration by parts, it is easy to see that $\la{\mathbb
L}\nabla^2\psi\ra$ is completely determined by the boundary data
$\{u_3,u_{3,n},({\mathbb L}\nabla^2\psi){\bf n}\cdot{\bf n},\
\text{div}({\mathbb L}\nabla^2\psi)\cdot{\bf n}+(({\mathbb
L}\nabla^2\psi){\bf n}\cdot{\bf t})_{,t}\}$.

On the other hand, simple integration by parts implies
\begin{equation}\label{1}
\int_{\Omega}\text{det}\nabla^2\psi=\int_{\partial\Omega}\psi_{,1}\psi_{,22}n_1-\psi_{,1}\psi_{,12}n_2=\int_{\partial\Omega}\psi_{,2}\psi_{,11}n_2-\psi_{,2}\psi_{,12}n_1.
\end{equation}
Recall that ${\bf t}=-{\bf R}_{\perp}{\bf n}$. Then we obtain
from \eqref{1} that
\begin{equation}\label{2}
\la\text{det}\nabla^2\psi\ra=\frac{1}{2|\Omega|}\int_{\partial\Omega}(\nabla^2\psi){\bf R}_{\perp}{\bf
n}\cdot {\bf R}_{\perp}\nabla
\psi=-\frac{1}{2|\Omega|}\int_{\partial\Omega}(\nabla^2\psi){\bf t}\cdot
{\bf R}_{\perp}\nabla \psi.
\end{equation}
In view of the expression of $\nabla \psi$ on $\partial\Omega$,
we have that
\begin{equation}\label{3}
{\bf R}_{\perp}\nabla \psi=\psi_{,n}{\bf R}_{\perp}{\bf n}+\psi_{,t}{\bf R}_{\perp}{\bf
t}=-\psi_{,n}{\bf t}+\psi_{,t}{\bf n}.
\end{equation}
Substituting \eqref{3} into \eqref{2} gives
\begin{equation}\label{4}
\la\text{det}\nabla^2\psi\ra=\frac{1}{2|\Omega|}\int_{\partial\Omega}\psi_{,n}(\nabla^2\psi){\bf
t}\cdot {\bf
t}-\frac{1}{2|\Omega|}\int_{\partial\Omega}\psi_{,t}(\nabla^2\psi){\bf
t}\cdot{\bf n}.
\end{equation}
Easy computations show that
$$(\nabla^2\psi){\bf t}\cdot {\bf t}=\psi_{,tt}\quad\text{and}\quad(\nabla^2\psi){\bf
t}\cdot{\bf n}=\psi_{,nt}=(\psi_{,n})_{,t}$$ and \eqref{4} is equivalent
to
\begin{eqnarray*}
\la\text{det}\nabla^2\psi\ra&=&\frac{1}{2|\Omega|}\int_{\partial\Omega}\psi_{,n}\psi_{,tt}-\frac{1}{2|\Omega|}\int_{\partial\Omega}\psi_{,t}(\psi_{,n})_{,t}\\
&=&\frac{1}{|\Omega|}\int_{\partial\Omega}\psi_{,n}\psi_{,tt}=-\frac{1}{|\Omega|}\int_{\partial\Omega}(\psi_{,n})_{,t}\psi_{,t}.
\end{eqnarray*}
In other words, $\la\text{det}\nabla^2\psi\ra$ is completely determined
by the Dirichlet data
$\{\psi,\psi_{,n}\}$.

Our next observation is that
\[
\Big{\la}\frac 12\div\div\left({\mathbb
R}(\nabla\theta\otimes\nabla u_3)+{\mathbb
R}(\nabla\theta\otimes\nabla u_3)^T\right)\psi+\div({\mathbb
R}\nabla^2\psi\nabla\theta)u_3\Big{\ra}
\]
is a null-Lagrangian as well. Indeed, using the integration by parts, one can easily check that
\begin{eqnarray}\label{18}
&&-\Big{\la}\frac 12\div\div\left({\mathbb
R}(\nabla\theta\otimes\nabla u_3)+{\mathbb
R}(\nabla\theta\otimes\nabla u_3)^T\right)\psi+\div({\mathbb R}\nabla^2\psi\nabla\theta)u_3\Big{\ra}\notag\\
&=&-\Big{\la}\frac 12\left({\mathbb R}(\nabla\theta\otimes\nabla
u_3)+{\mathbb R}(\nabla\theta\otimes\nabla
u_3)^T\right)\cdot\nabla^2\psi-{\mathbb
R}\nabla^2\psi\nabla\theta\cdot\nabla
u_3\Big{\ra}\notag\\
&&+\frac{1}{|\Omega|}\int_{\partial\Omega}\left\{\left(-\frac 12\div\left({\mathbb
R}(\nabla\theta\otimes\nabla u_3)+{\mathbb
R}(\nabla\theta\otimes\nabla u_3)^T\right)\cdot{\bf n}\right.\right.\notag\\
&&\hspace{10mm}\left.-\left(\frac 12\left({\mathbb
R}(\nabla\theta\otimes\nabla u_3)+{\mathbb
R}(\nabla\theta\otimes\nabla u_3)^T\right){\bf n}\cdot{\bf t}\right)_{,t}\right)\psi\notag\\
&&\hspace{10mm}+\left.\left(\frac 12\left({\mathbb
R}(\nabla\theta\otimes\nabla u_3)+{\mathbb
R}(\nabla\theta\otimes\nabla u_3)^T\right){\bf n}\cdot{\bf n}\right)\psi_{,n}-{\mathbb R}\nabla^2\psi\nabla\theta\cdot{\bf n}u_3\right\}\notag\\
&=&-\big{\la}(\nabla\theta\otimes\nabla u_3)\cdot{\mathbb
R}\nabla^2\psi-{\mathbb R}\nabla^2\psi\nabla\theta\cdot\nabla
u_3\big{\ra}+B\notag\\
&=&B,
\end{eqnarray}
where $B$ is the boundary integral term given above.  Thus, $B$ is determined by $\psi,\psi_{,n},u_3,u_{3,n}$.

\section{Bounds on the volume fraction}

Now we are ready to derive bounds on the volume fraction of 2-phase shallow shell in \eqref{shell} by one set of Cauchy data:
\begin{equation}\label{data2}
\{({\bf u}';{u}_3, u_{3,n}), ({\bf s}^{\theta}{\bf n};
(\div\ {\bf m}-{\bf s}^{\theta}\nabla\theta){\bf n}+({\bf m}{\bf n}\cdot {\bf t})_{,t}, {\bf m}{\bf n}\cdot {\bf n})\}.
\end{equation}
From the above result, we can see that this problem is reduced to the same problem for \eqref{shell2} with corresponding Cauchy data, which is completely determined by \eqref{data2}. Moreover, the boundary values of $\{({u}_3, u_{3,n});(\div({\mathbb M}\nabla^2u_3)\cdot{\bf n}+(({\mathbb M}\nabla^2u_3){\bf n}\cdot {\bf t})_{,t}, {\mathbb M}\nabla^2u_3{\bf n}\cdot {\bf n})\}$ and $\{(\psi, \psi_{,n});(\div({\mathbb L}\nabla^2\psi)\cdot{\bf n}+(({\mathbb L}\nabla^2\psi){\bf n}\cdot {\bf t})_{,t}, {\mathbb L}\nabla^2\psi{\bf n}\cdot {\bf n})\}$  are also determined by \eqref{data2}. Therefore, using \eqref{shell2}, \eqref{18}, and integration by parts, we obtain that
$\la{\mathbb L}\nabla^2\psi\cdot\nabla^2\psi\ra+\la{\mathbb M}\nabla^2u_3\cdot\nabla^2u_3\ra$
is a null-Lagrangian. Denote $\la{\mathbb L}\nabla^2\psi\ra=b_0, \la{\mathbb M}\nabla^2u_3\ra=\tilde b_0,$
\begin{equation*}
{\mathcal A}=\{\la{\mathbb L}\nabla^2w\ra=b_0;\ w=\psi,\ w_{,n}=\psi_{,n}\ \text{on}\ \partial\Omega\},
\end{equation*}
and
\begin{equation*}
{\mathcal B}=\{\la{\mathbb M}\nabla^2v\ra=\tilde b_0;\ v=u_3,\ v_{,n}=u_{3,n}\ \text{on}\ \partial\Omega\}.
\end{equation*}
It is obvious that
\begin{eqnarray}\label{3.2}
e_0&:=&\la{\mathbb L}\nabla^2\psi\cdot\nabla^2\psi\ra+\la{\mathbb M}\nabla^2u_3\cdot\nabla^2u_3\ra\notag\\
&\ge&\underset{\underline{\psi}\in{\mathcal A}}{\min}\la{\mathbb L}\nabla^2{\underline{\psi}}\cdot\nabla^2{\underline{\psi}}\ra+\underset{\underline{u}_3\in{\mathcal B}}{\min}\la{\mathbb M}\nabla^2{\underline{u}}_3\cdot\nabla^2{\underline{u}}_3\ra.
\end{eqnarray}

It should be noted that $\underset{\underline{\psi}\in{\mathcal A}}{\min}\la{\mathbb L}\nabla^2{\underline{\psi}}\cdot\nabla^2{\underline{\psi}}\ra$ (resp. $\underset{\underline{u}_3\in{\mathcal B}}{\min}\la{\mathbb M}\nabla^2{\underline{u}}_3\cdot\nabla^2{\underline{u}}_3\ra$) does not necessarily occur at the solution of the corresponding thin plate equation $\div\div({\mathbb L}\nabla^2w)=0$ with Dirichlet data $w=\psi$ and $w_{,n}=\psi_{,n}$ (resp.  $\div\div({\mathbb M}\nabla^2v)=0$ with Dirichlet data $v=u_3$ and $v_{,n}=u_{3,n}$) since we do not know whether $\la{\mathbb L}\nabla^2w\ra=b_0$ (resp. $\la{\mathbb M}\nabla^2v\ra=\tilde b_0$) holds. To be precise,  the boundary value problem 
\begin{equation*}
\begin{cases}
\div\div({\mathbb L}\nabla^2w)=0\ \text{in}\ \Omega,\\
w=\psi,\ w_{,n}=\psi_{,n}\ \text{on}\ \partial\Omega
\end{cases}
\end{equation*}
 has a unique solution and thus the corresponding Neumann data $\{\div({\mathbb L}\nabla^2w)\cdot{\bf n}+(({\mathbb L}\nabla^2w){\bf n}\cdot {\bf t})_{,t}, ({\mathbb L}\nabla^2w){\bf n}\cdot {\bf n}\}$ are completely determined by $\psi,\psi_{,n}$. Recall that $\la{\mathbb L}\nabla^2w\ra$ depends on the Neumann data.  It is not known whether $\la{\mathbb L}\nabla^2w\ra=\la{\mathbb L}\nabla^2\psi\ra(=b_0)$. However, we have shown that for $\underline{\psi}\in{\mathcal A}$ (resp. $\underline{u}_3\in{\mathcal B}$) we can determine $\la\nabla^2\underline{\psi}\ra:=a_0$ (resp. $\la\nabla^2\underline{u}_3\ra:=\tilde a_0$) and $\la\text{det}\nabla^2\underline{\psi}\ra:=c_0$ (resp. $\la\text{det}\nabla^2\underline{u}_3\ra:=\tilde c_0$). Therefore, it follows from \eqref{3.2} that
\begin{equation}\label{min0}
\begin{split}
e_0-2\zeta c_0-2\tilde\zeta\tilde c_0&\ge\underset{\underline{\psi}\in{\mathcal A}}{\min}\la{\mathbb L}\nabla^2{\underline{\psi}}\cdot\nabla^2{\underline{\psi}}\ra-2\zeta c_0+\underset{\underline{u}_3\in{\mathcal B}}{\min}\la{\mathbb M}\nabla^2{\underline{u}}_3\cdot\nabla^2{\underline{u}}_3\ra-2\tilde\zeta\tilde c_0\\
&\ge\underset{\la\text{det}\underline{A}\ra=c_0}{\underset{\la{\mathbb L}\underline{A}\ra=b_0}{\underset{\la\underline{A}\ra=a_0}{\min}}}\la \underline{A}\cdot {\mathbb L}
\underline{A}\ra-2\zeta c_0+\underset{\la\text{det}\underline{B}\ra=\tilde c_0}{\underset{\la{\mathbb M}\underline{B}\ra=\tilde b_0}{\underset{\la\underline{B}\ra=\tilde a_0}{\min}}}\la \underline{B}\cdot {\mathbb M}
\underline{B}\ra-2\tilde\zeta \tilde c_0,
\end{split}
\end{equation}
where $\zeta$ and $\tilde\zeta$ are two parameters which will be specified later. It is clear that the two minimization problems on the right hand side of \eqref{min0} can be treated separately. 

We first consider 
\begin{equation}\label{min00}
\underset{\la\text{det}\underline{A}\ra=c_0}{\underset{\la{\mathbb L}\underline{A}\ra=b_0}{\underset{\la\underline{A}\ra=a_0}{\min}}}\la \underline{A}\cdot {\mathbb L}
\underline{A}\ra-2\zeta c_0.
\end{equation}
Introducing the basis
$${\mathcal P}=\left\{\frac{1}{\sqrt{2}}\begin{pmatrix}1&0\\0&1\end{pmatrix},\frac{1}{\sqrt{2}}\begin{pmatrix}1&0\\0&-1\end{pmatrix},\frac{1}{\sqrt{2}}\begin{pmatrix}0&1\\1&0\end{pmatrix}\right\},$$
then any symmetric matrix
$$\underline{A}=\begin{pmatrix}a_{11}&a_{12}\\a_{12}&a_{22}\end{pmatrix}$$ is
represented by a three-dimensional vector
${\bf v}=\frac{1}{\sqrt{2}}(a_{11}+a_{22},a_{11}-a_{22},2a_{12})^T$ and
the determinant of $\underline{A}$ is written as
$$\text{det}\underline{A}=\frac 12{\bf v}\cdot T{\bf v}$$ with
$$T=\begin{pmatrix}1&0&0\\0&-1&0\\0&0&-1\end{pmatrix}.$$ Similarly, in the basis ${\mathcal P}$, the 4th tensor ${\mathbb L}$ is expressed by the $3\times 3$ matrix
$$L=\begin{pmatrix}\alpha&0&0\\0&\beta&0\\0&0&\beta\end{pmatrix}$$ with
\[
\alpha=\frac{\lambda+2\mu}{4\mu(2\mu+3\lambda)}\quad\text{and}\quad\beta=\frac{1}{4\mu}.
\]
Assume the strong convexity condition is satisfied, i.e., $\mu>0$ and $2\mu+3\lambda>0$. We can see that both $\alpha$ and $\beta$ are positive. We denote $\alpha_1,\beta_1$ and $\alpha_2,\beta_2$ the values in phase 1 and phase 2, respectively. Define $\alpha_{\ast}=\min\{\alpha_1,\alpha_2\}$ and $\beta_{\ast}=\min\{\beta_1,\beta_2\}$.  Introduce the translated tensor
\[
S=L-\zeta T
\]
for $\zeta\in(-\beta_{\ast},\alpha_{\ast})$, then $S$ is also positive. Denote $g_0=b_0-\zeta Ta_0$. Then the minimization problem in \eqref{min00} is reduced to
\begin{equation}\label{min2}
\underset{\la\text{det}\underline{A}\ra=c_0}{\underset{\la{\mathbb L}\underline{A}\ra=b_0}{\underset{\la\underline{A}\ra=a_0}{\min}}}\la \underline{A}\cdot {\mathbb L}
\underline{A}\ra-2\zeta c_0\ge\underset{\underset{\la S\underline{A}\ra=g_0}{{\la
\underline{A}\ra=a_0}}}{\min}\la \underline{A}\cdot S
\underline{A}\ra.
\end{equation}

The minimization problem on the right hand side of \eqref{min2} is
exactly the one studied in Milton and Nguyen \cite{ml} (see (3.8)
there). Using their result, if $\hat{A}$ is the minimizer, then
\begin{equation}\label{ml1}
\la\hat{A}\cdot S\hat{A}\ra=g_0\cdot
a_0+\frac{(S_1-S_2)^{-2}}{f_1f_2}[g_0-\la S\ra
a_0]\cdot[(f_2S_1+f_1S_2)g_0-S_1S_2a_0].
\end{equation}
Applying the same arguments to 
\begin{equation}\label{ml10}
\underset{\la\text{det}\underline{B}\ra=\tilde c_0}{\underset{\la{\mathbb M}\underline{B}\ra=\tilde b_0}{\underset{\la\underline{B}\ra=\tilde a_0}{\min}}}\la \underline{B}\cdot {\mathbb M}
\underline{B}\ra-2\tilde\zeta \tilde c_0\ge\underset{\underset{\la\tilde S\underline{B}\ra=\tilde g_0}{{\la
\underline{B}\ra=\tilde a_0}}}{\min}\la \underline{B}\cdot \tilde S
\underline{B}\ra,
\end{equation}
we obtain that if $\hat{B}$ is the minimizer then
\begin{equation}\label{ml2}
\la\hat{B}\cdot \tilde S\hat{B}\ra=\tilde g_0\cdot
\tilde a_0+\frac{(\tilde S_1-\tilde S_2)^{-2}}{f_1f_2}[\tilde g_0-\la \tilde S\ra
\tilde a_0]\cdot[(f_2\tilde S_1+f_1\tilde S_2)\tilde g_0-\tilde S_1\tilde S_2\tilde a_0],
\end{equation}
where
\begin{equation*}
\left\{\begin{split}
&\tilde S=M-\tilde\zeta T,\\
&M=\begin{pmatrix}\tilde\alpha&0&0\\0&\tilde\beta&0\\0&0&\tilde\beta\end{pmatrix},\ \tilde\alpha=\frac{4\mu(3\lambda+2\mu)}{3(\lambda+2\mu)},\ \tilde\beta=\frac{4\mu}{3},\\
&\tilde\zeta\in(-\tilde\beta_{\ast},\tilde\alpha_{\ast}),\ \tilde\alpha_{\ast}=\min\{\tilde\alpha_1,\tilde\alpha_2\},\ \tilde\beta_{\ast}=\min\{\tilde\beta_1,\tilde\beta_2\},\\
&\tilde g_0=\tilde b_0-\tilde\zeta T\tilde a_0.
\end{split}\right.
\end{equation*}

Combining \eqref{min0}, \eqref{min2}, \eqref{ml1}, \eqref{ml10}, and \eqref{ml2} gives
\begin{equation}\label{min3}
\begin{split}
&\quad\ e_0-2\zeta c_0-b_0\cdot a_0+2\zeta\text{det}a_0-2\tilde\zeta \tilde c_0-\tilde b_0\cdot \tilde a_0+2\tilde\zeta\text{det}\tilde a_0\\
&\ge\frac{(S_1-S_2)^{-2}}{f_1f_2}[g_0-\la S\ra
a_0]\cdot[(f_2S_1+f_1S_2)g_0-S_1S_2a_0]\\
&\quad\ +\frac{(\tilde S_1-\tilde S_2)^{-2}}{f_1f_2}[\tilde g_0-\la \tilde S\ra
\tilde a_0]\cdot[(f_2\tilde S_1+f_1\tilde S_2)\tilde g_0-\tilde S_1\tilde S_2\tilde a_0].
\end{split}
\end{equation}
Expanding the first term on the right hand side of \eqref{min3} reveals
\begin{eqnarray}\label{min4}
&&\frac{(S_1-S_2)^{-2}}{f_1f_2}[g_0-\la S\ra
a_0]\cdot[(f_2S_1+f_1S_2)g_0-S_1S_2a_0]\notag\\
&=&\frac{(S_1-S_2)^{-2}}{f_1f_2}[(b_0-\zeta Ta_0)-\la L\ra a_0+\zeta
Ta_0]\notag\\
&&\cdot[(f_2L_1+f_1L_2-\zeta T)(b_0-\zeta Ta_0)-(L_1L_2-\zeta
TL_2-\zeta L_1T+\zeta^2 T^2)a_0]\notag\\
&=&\frac{(S_1-S_2)^{-2}}{f_1f_2}[b_0-\la L\ra a_0]\cdot[(f_2L_1+f_1L_2)b_0-\zeta Tb_0-\zeta(f_2L_1+f_1L_2)Ta_0\notag\\
&&-(L_1L_2-\zeta TL_2-\zeta L_1T)a_0].
\end{eqnarray}
Therefore, this term is linear in $\zeta$ with leading coefficient
\begin{eqnarray}\label{min5}
&&\frac{(S_1-S_2)^{-2}}{f_1f_2}[b_0-\la L\ra
a_0]\cdot[(TL_2+L_1T)a_0-Tb_0-(f_2L_1+f_1L_2)Ta_0]\notag\\
&=&\frac{(S_1-S_2)^{-2}}{f_1f_2}[b_0-\la L\ra
a_0]\cdot[f_1L_1Ta_0+f_2L_2Ta_0-Tb_0]\notag\\
&=&\frac{(S_1-S_2)^{-2}}{f_1f_2}[b_0-\la L\ra a_0]\cdot[\la L\ra
Ta_0-Tb_0].
\end{eqnarray}
Note that $L_2$ and $T$ commute. Likewise, we can expand the second term on the right hand side of \eqref{min3} and obtain that it is linear in $\tilde\zeta$ with leading coefficient
\begin{equation}\label{min51}
\frac{(\tilde S_1-\tilde S_2)^{-2}}{f_1f_2}[\tilde b_0-\la M\ra \tilde a_0]\cdot[\la M\ra
T\tilde a_0-T\tilde b_0].
\end{equation}
Putting \eqref{min3}-\eqref{min51} together, we thus prove that
\begin{theorem}\label{bounds}
The following bounds hold:
\begin{equation}\label{min6}
\begin{split}
&e_0-b_0\cdot a_0-\tilde b_0\cdot\tilde a_0\\
&\ge\zeta\left\{\frac{(S_1-S_2)^{-2}}{f_1f_2}[b_0-\la L\ra
a_0]\cdot[\la L\ra Ta_0-Tb_0]+2c_0-2\text{\rm det}a_0\right\}\\
&+\frac{(S_1-S_2)^{-2}}{f_1f_2}[b_0-\la L\ra
a_0]\cdot[(f_2L_1+f_1L_2)b_0-L_1L_2a_0]\\
&+\tilde\zeta\left\{\frac{(\tilde S_1-\tilde S_2)^{-2}}{f_1f_2}[\tilde b_0-\la M\ra
\tilde a_0]\cdot[\la M\ra T\tilde a_0-T\tilde b_0]+2\tilde c_0-2\text{\rm det}\tilde a_0\right\}\\
&+\frac{(\tilde S_1-\tilde S_2)^{-2}}{f_1f_2}[\tilde b_0-\la M\ra
\tilde a_0]\cdot[(f_2M_1+f_1M_2)\tilde b_0-M_1M_2\tilde a_0]
\end{split}
\end{equation}
for all $\zeta\in(-\beta_{\ast},\alpha_{\ast})$ and $\tilde\zeta\in(-\tilde\beta_{\ast},\tilde\alpha_{\ast})$. 
\end{theorem}

We remark that \eqref{min6} is also valid at the extremal
points $(\zeta,\tilde\zeta)=(-\beta_{\ast},-\tilde\beta_{\ast})$, $(-\beta_{\ast},\tilde\alpha_{\ast})$, $(\alpha_{\ast},-\tilde\beta_{\ast})$, or $(\alpha_{\ast},\tilde\alpha_{\ast})$. Therefore, by substituting the extreme points into \eqref{min6}, we obtain the tightest bounds on the volume fraction $f_1$ (or $f_2$) given by these inequalities. Moreover, when the shell is flat ($\theta=0$) and the plate measurement is either under pure stretching ($u_3=0$) or pure bending ($u_1=u_2=0$) then the bounds derived here reduce to those derived in \cite{ml}. Note that the pure bending case is equivalent to the pure stretching case as observed in \cite{kmw}. 

\section*{Acknowledgements}
HK was partially supported by National Research Foundation of Korea
through NRF grants No. 2009-0085987 and 2010-0017532. GWM was
partially supported by the National Science Foundation of the USA
through grant DMS-0707978. JNW was partially supported by the
National Science Council of Taiwan through grants 100-2628-M-002-017
and 99-2115-M-002-006-MY3.


\begin{thebibliography}{999}

\bibitem{cm}
P.G. Ciarlet and B. Miara, \emph{Justification of the
two-dimensional equations of a linearly elastic shallow shell},
Comm. Pure Appl. Math., \textbf{45} (1992), 327-360.

\bibitem{dlw}
M. Di Cristo, C.-L. Lin, and J.-N. Wang, \emph{Quantitative uniqueness estimates for the shallow shell system and their application to an inverse problem}, to appear in Ann. Sc. Norm. Super. Pisa Cl. Sci.

\bibitem{kkm}
H. Kang, E. Kim, and G.W. Milton, \emph{Sharp bounds on the volume fractions of two materials in a two-dimensional body from electrical boundary measurements: the translation method}, Cal. Var. Partial Diff. Equations, (2011).

\bibitem{kmw}
H. Kang, G. Milton, and J.N. Wang, \emph{Equivalence of inverse problems for the plane elasticity and the thin plate with finite measurements and its applications}, arXiv:1203.3833v1 [math.AP].

\bibitem{lc1}
K.A. Lurie and A.V. Cherkaev, \emph{Accurate estimates of the conductivity of mixtures formed of two materials in a given proportion (two-dimensional problem)}, Doklady Akademii Nauk SSSR, \textbf{264} (1982), 1128-1130. English translation in Soviet Phys. Dokl., \textbf{27} (1982), 461-462.

\bibitem{lc2}
K.A. Lurie and A.V. Cherkaev, \emph{Exact estimates of conductivity of composites formed by two isotropically conducting media taken in prescribed proportion}, Proc. Royal Soc. Edinburgh. Section A, Math. Phys. Sciences, \textbf{99} (1984), 71-87.

\bibitem{ml}
G.W. Milton and L.H. Nguyen, \emph{Bounds on the volume fraction of
2-phase, 2-dimensional elastic bodies and on (stress, strain) pairs
in composite}, Comptes rendus - M\'ecanique, Reference: CRAS2B3066, DOI information: 10.1016/j.crme.2012.02.002.

\bibitem{mbook}
G.W. Milton, \emph{The Theory of Composites}, Cambridge Monographs
on Applied and Computational Mathematics, Cambridge University
Press, 2002.

\bibitem{mt1}
F. Murat and L. Tartar, \emph{Calcul des variations et
homog\'en\'isation, (French) {\rm [}Calculus of variation and
homogenization{\rm ]}}, In Les m\'ethodes de l'homog\'en\'eisation:
th\'eorie et applications en physique, volume 57 of Collection de la
Direction des \'etudes et recherches d'Electricit\'e de France,
pages (1985)319-369, Paris, Eyrolles. English translation in Topics
in the Mathematical Modelling of Composite Materials, 139-173, ed.
by A. Cherkaev and R. Kohn, ISBN 0-8176-3662-5.


\bibitem{mt2}
L. Tartar, \emph{Estimation de coefficients homog\'en\'eis\'es.
(French) {\rm [}Estimation of homogenization coefficients{\rm ]}}.
In R. Glowinski and J.-L. Lions, editors, Comp. Meth. Appl. Sci.
Eng.: Third International Symposium, Versailles, France, December
5-9 (1977), volume 704 of Lecture Notes in Mathematics, pages
(1979) 364-373, Berlin, Springer-Verlag. English translation in
Topics in the Mathematical Modelling of Composite Materials, 9-20,
ed. by A. Cherkaev and R. Kohn. ISBN 0-8176-3662-5.


\bibitem{mt3}
L. Tartar, \emph{Estimations fines des coefficients
homog\'en\'eis\'es, (French) {\rm [}Fine estimations of homogenized
coefficients{\rm ]}}, In P. Kr\'ee, editor, Ennio de Giorgi
Colloquium: Papers Presented at a Colloquium Held at the H.
Poincar\'e Institute in November 1983, volume 125 of Pitman Res.
Notes in Mathematics, (1985) 168-187, London. Pitman Publishing
Ltd.

\end{thebibliography}
\end{document}